\documentclass[12pt, amstex, amssymb]{article}
\input{diagrams.sty}
\newarrow{Rel}{-}{-}{+}{-}{>}
\newarrow{Hook}{C}{-}{-}{-}{>}
\title{Elements of a theory of algebraic theories}
\author{J. M. E. Hyland}
\newtheorem{theorem}{Theorem}[section]
\newtheorem{proposition}[theorem]{Proposition}

\newtheorem{definition}[theorem]{Definition}
\def\slashedrightarrow{\relbar\joinrel\mapstochar\joinrel\rightarrow}
\newcommand{\eps}{\varepsilon}

\newcommand{\relto}{\slashedrightarrow}
\newcommand{\opcomm}{\odot}
\newcommand{\Mult}[1]{#1{\text{-Mult}}}
\newcommand{\Alg}[1]{#1{\text{-Alg}}}
\newcommand{\PsAlg}[1]{#1{\text{-PsAlg}}}
\newcommand{\Kleisli}[1]{{\text{Kl}}(#1)}
\newcommand{\cell}{\Longrightarrow}
\usepackage{amsmath, amssymb}
\usepackage[all,cmtip]{xy}

\begin{document}

\maketitle

\begin{abstract}
Kleisli bicategories are a natural environment
in which the combinatorics involved in various notions of 
algebraic theory can be handled in a uniform way.
The setting allows a clear account of comparisons between such
notions. Algebraic theories, symmetric operads and
nonsymmetric operads are treated as examples.
\end{abstract}

\section{Introduction}
This paper has its genesis in Glynn Winskel's use of
presheaf categories and profunctors in the foundations of 
concurrency. His basic theory is laid out in \cite{CW05}
with Cattani, and particular cases of Kleisli bicategories
appear there. A preordered set version, providing a model
for linear logic, is already in \cite{WZ04}.
Kleisli bicategories are both a rich source of models and
a context in which to understand subtle theory.
Their value was recognised by a group of us in Cambridge and 
we set about preparing an exposition \cite{FGHW} of the general theory.
Around the same time John Power realised the significance
of the key pseudo-distributivities
in connection with extensions of Edinburgh
work \cite{FPT99} on variable binding. The paper \cite{CHP03}
shows the common interest and in Edinburgh a thesis \cite{Tanaka05} and 
papers (for example \cite{TP06}) quickly followed.
By contrast the Cambridge
exposition remains unfinished, and
there is
just one paper \cite{FGHW08} which gives some
sense of our preferred approach. That is my
fault and I have written this paper for Glynn Winskel by way of
apology. It is not intended as a substitute for the unfinished paper.
Rather it sketches applications to 
algebraic theories and operads, which I have presented
in talks over the years. In developing the ideas,
I have profited from discussions on with Richard Garner and John Power.
Recently Garner and I have made progress
on coalgebraic aspects, and a substantial theory is emerging.
Here I focus on just one strand of ideas, and leave details
and the wider perspective 
for other occasions.

The paper is organised as follows.
In section 2, I describe and give elementary
properties of the basic construction, 
that of the Kleisli bicategory $\Kleisli{P}$ of
a Kleisli structure $P$. This is in my view a good
way to understand the bicategories with which we 
shall be concerned, and I explain how even the
basic bicategory ${\bf Prof}$ of profunctors,
which underlies the whole paper, can be considered
from this perspective.
Section 3 is concerned with distributive laws,
composite Kleisli structures and features of the
corresponding Kleisli categories. The idea of 
a pseudo-distributive law at roughly the level 
we need is old, see \cite{Kelly74}.
Perhaps because there were no compelling applications,
the details were not worked through. The first 
complete account seems to be \cite{Marm99}. 
The use here of Kleisli structures creates
a new focus but there are no great surprises.

In section 4, I describe my approach to
algebraic theories as monads in Kleisli bicategories.
A concrete categorical treatment of essentially the same point of view
is in \cite{Cur12}. The value of an abstract treatment
becomes more apparent with very recent work but even at the level of
this paper I hope readers
will appreciate the smooth treatment of categories
generated by theories. 
In the final Section 5, I use the general setting to 
give a treatment of comparisons between notions of algebraic theory.
I hope inter alia to encourage sensitivity to
some subtleties in the notions of symmetric
and non-symmetric operad.

I make a small remark about notation. I have decided to use
$P$ both when describing a general Kleisli structure
and for the presheaf Kleisli structure. I hope to avoid confusion
by the following convention. In the general case I use
lower case letter for the objects of a bicategory.
In the special case when the objects are themselves
categories I use upper case letters.

These introductory remarks make clear that I have discussed
material with many people. But on this occasion
I thank in
particular Glynn Winskel. We go back a long way
and it has been continually stimulating to discuss logic
and computer science with him over many years. I very much
appreciate his intellectual honesty and openness, and his
talent for grounding abstract mathematics in the modelling of 
computational phenomena. This paper derives from work of
his and he was the first person with whom I discussed
Kleisli bicategories. I dedicate the paper to him
with affection and best wishes for the future.

\section{Kleisli Bicategories}

The Kleisli formulation of a monad on a category $\mathbb C$,
given in \cite{Manes76},
has not played a prominent role in mathematics. 
but it is familiar 
in the programming language community from   
the computational $\lambda$-calculus \cite{Moggi89} and 
premonoidal categories \cite{PR97}.
The Kleisli presentation gives for each object $c$ in a category
${\mathbb C}$ a unit $c \to Tc$ in
$\mathbb C$ and for each $f : c \to Td$ in ${\mathbb C}$ a lift
$f^{\sharp} : Tc \to Td$. This data is required to satisfy
evident equations. The virtue of the formulation
is that it makes trivial the definition and basic
properties of the Kleisli category of a monad.
Now one can still define a Kleisli
category when structure is only given for some subcollection
of the objects of $\mathbb C$ and generalisations of this
kind have been identified, for example in \cite{ACU12}. 
The phrase relative monad
is in use, but for the cases considered 
here I prefer to say restricted.

\subsection{Kleisli structures}
A Kleisli structure is a $2$-dimensional version
of a restricted monad.
The starting point is a 
bicategory $\cal K$
equipped with a sub-bicategory 
${\cal A} \hookrightarrow {\cal K}$
\begin{definition} \label{def:kles}  
A {\em Kleisli structure} $P$ on ${\cal A} \hookrightarrow {\cal K}$
is the following.
\begin{itemize}
\item A choice for each object $a \in {\cal A}$, of
an arrow $y_a : a \rightarrow Pa$ in $\cal K$.
\item For each pair $a, b \in {\cal A}$ of objects, a functor
\[
 {\cal K}(a, Pb)  \longrightarrow  {\cal K}(Pa, Pb)  \qquad
  f \longmapsto      f^\sharp 
\] 
\item Families of invertible 2-cells
\[
\eta_f : f \rightarrow f^\sharp . y_a \qquad
 \kappa_a : (y_a)^\sharp \rightarrow 1_{Pa} \qquad
\kappa_{g, f} : (g^\sharp . f)^\sharp \rightarrow g^\sharp . f^\sharp
\]
natural in $f : a \rightarrow Pb$ and $g: b \rightarrow Pc$
as appropriate,
and subject to unit and pentagon coherence conditions.
\end{itemize}
\end{definition}
It is clear from the data that $P$ can be given the 
structure of a pseudo-functor
$P: {\cal A} \to {\cal K}$. For 
$f: a \to b$ in $\cal A$, set 
$Pf = (y_bf)^{\sharp} : Pa \to Pb$. The 2-cell structure
and its coherence are routine. Then
$y_a : a \to Pa$ can be given the structure of a
transformation. For $a \stackrel{f}{\to} b$ the structure 
$2$-cells are
$y_b . f \stackrel{\eta_{y_b. f}}{\longrightarrow}
(y_b . f)^\sharp . y_a = Pf . y_a $.

\subsection{The Kleisli bicategory}\label{bicatsec}
Given a Kleisli 
structure $P$ on ${\cal A} \hookrightarrow {\cal K}$
we define its {\em Kleisli bicategory} $\Kleisli{P}$ as follows.
The objects of $\Kleisli{P}$ are the objects of $\cal A$. 
For objects $a, b$, set $\Kleisli{P}(a,b) = {\cal K}(a,Pb)$.
The identities of $\Kleisli{P}$ are the $y_a : a \rightarrow Pa$, 
from the Kleisli structure. 
The Kleisli composition 
of $f : a \rightarrow Pb$ and 
$g : b \rightarrow Pc$, is $g \cdot f= g^\sharp. f : a \to Pc$. 
This extends to 2-cells so we have 
composition functors. To obtain a bicategory, it remains to define 
coherent unit and associativity isomorphisms 
$\lambda_f : y_b \cdot f \rightarrow f $, 
$\rho_f :  f \rightarrow f \cdot y_a $ and
$\alpha_{h, g, f}:   (h \cdot g) \cdot f \rightarrow h \cdot (g \cdot f)$.
The unit isomorphisms $\lambda_f$ and $\rho_f$ are given
by $(y_b)^{\sharp}f \stackrel{\kappa_b f}{\longrightarrow} 1_{Pb}f \cong f$
and $f \stackrel{\eta_f}{\longrightarrow} f^{\sharp}y_a$
respectively, while the associativity $\alpha_{h, g, f}$
is
$(h^{\sharp}g)^{\sharp} f \stackrel{\kappa_{h,g}f}{\longrightarrow} 
(h^{\sharp}g^{\sharp}f \cong h^{\sharp}(g^{\sharp}f)$.
The coherence axioms follow directly from the coherence conditions
of the Kleisli structure.
\begin{theorem}\label{kleislibicat} Let $P$ be a Kleisli 
structure on ${\cal A} \hookrightarrow {\cal K}$. Then $\Kleisli{P}$ is 
a bicategory.
\end{theorem}
For simplicity in what follows, I shall adopt standard notation
and often write $a \relto b$ instead of $a \to PB$ for maps
in $\Kleisli{P}(a,b)$.

In traditional category theory, the Kleisli construction is one universal
way to associate an adjunction with a monad. In the $2$-dimensional
setting of Kleisli structures we get a restricted (pseudo)adjunction
as follows. There is a `forgetful' pseudofunctor $\Kleisli{P} \to {\cal K}$
taking $f: a \to Pb$ in $\Kleisli{P}$ to $f^{\sharp}: Pa \to Pb$ in 
$\cal K$. And there is a pseudofunctor $F: {\cal A} \to \Kleisli{P}$,
taking $f: a \to b$ in $\cal A$ to $y_b.f : a \to Pb$ considered as
a map from $a$ to $b$ in $\Kleisli{P}$. I omit the $2$-dimensional 
structure which is routine, but note that the fact that $F$ is a 
restricted left pseudoadjoint is immediate from the identification
$\Kleisli{P}(a,b) = {\cal K}(a,Pb)$.
  
\subsection{Presheaves and profunctors}
The basic example of a Kleisli structure is 
the {\em presheaf Kleisli structure} arising from the presheaf 
construction. This gives a Kleisli structure
on ${\bf Cat} \hookrightarrow {\bf CAT}$, the inclusion of the
$2$-category of small categories in that of locally small
categories. First recall my
notational convention: in general bicategories
objects are lower case; but now the objects are themselves
categories and I use upper case. 

Here are the details of the
presheaf Kleisli structure. For a small category $A$, take
$PA = [A^{op}, {\bf Sets}]$, the 
category of presheaves over $A$, with 
the usual Yoneda embedding
$y_A : A \rightarrow PA$. For $f : A \rightarrow PB$,
we have $f^\sharp :
PA \rightarrow PB$, a given choice
of left Kan extension of $f$ along the Yoneda
embedding. The $2$-dimensional structure arises
as part of the
story of why presheaves give a Kleisli
structure.

Let ${\bf COC}$ be the $2$-category of locally small
cocomplete categories, cocontinouous functors and
natural transformations. There is an evident forgetful
$2$-functor $U: {\bf COC} \to {\bf CAT}$. Consider
the diagram
\begin{diagram}
 & & {\bf COC}\\
 & \ruTo^{\hat{P}} & \dTo_{U}\\
{\bf Cat} & \rHook & {\bf CAT}
\end{diagram}
where $\hat{P}$ is the presheaf
construction thought of as a pseudo-functor ${\bf Cat} \to {\bf COC}$.
We have then for $A \in {\bf Cat}$ and ${\cal B} \in {\bf COC}$
an adjoint equivalence
\[
\xymatrix{
   {\bf CAT}(A, U{\cal B}) 
 \ar@<1.2ex>[rr]^(.52){(-)^\dagger}
 \ar@{}[rr]|{\bot}    
 & & \ar@<1.2ex>[ll]^{(-)  \, y_A} {\bf COC}(FA, {\cal B}) }
\]
where $(-)^{\dagger}$ is left Kan extension thought of as landing
in $\bf COC$. It is straightforward to extract Kleisli structure
from such a situation. We have $P = U\hat{P}$ as pseudo-functor
${\bf Cat} \to {\bf CAT}$ and the Yoneda 
%y_A : A \to PA$ is already in place. The Kleisli structure
left Kan extension is just $f^{\sharp} = U(f^{\dagger})$.
Now we can see how the $2$-cells in the Kleisli structure
arise. For $f: A \to PB = U\hat{P}B$, 
$\eta_f : f \to f^{\sharp}y_A = U(f^{\dagger})y_A$ is the unit
of the adjunction above. The $2$-cell
$\kappa_A : (y_A)^\sharp = U({y_A}^{\dagger}) \rightarrow 1_{PA}$ is
\begin{diagram} 
U( {y_A}^\dagger )\cong U\big((U(1_{FA}) y_A)^\dagger \big)&
\rTo{G(\varepsilon_{1_{FA}})} & U(1_{FA}) \cong 1_{PA}
\end{diagram}
using the counit of the adjunction.
The  $2$-cell 
$\kappa_{g, f} : (g^\sharp  f)^\sharp \to g^\sharp f^\sharp$
is the composite of the two lines displayed here.
\begin{diagram}
(g^\sharp f)^\sharp = G (G(g^\dagger) f)^\dagger &
\rTo^{G(G(g^\dagger) \eta_f)^\dagger} &  
G\big( G(g^\dagger) G(f^\dagger ) y_A \big)^\dagger
\cong G( G(g^\dagger f^\dagger)  y_A)^\dagger \, ;
\end{diagram}
\begin{diagram}
G( G(g^\dagger \cdot f^\dagger) \cdot y_A)^\dagger & 
\rTo^{G(\varepsilon_{g^\dagger \cdot f^\dagger})} &
G(g^\dagger \cdot f^\dagger) \cong G(g^\dagger) \cdot G(f^\dagger) \, .
\end{diagram}
With these definitions it is routine to check the Kleisli
structure coherence.

Consider now the Kleisli bicategory $\Kleisli{P}$ of the
presheaf Kleisli structure. It has objects the small categories,
and for $A, B \in {\bf Cat}$ we have
\[
\Kleisli{P} (A, B ) = {\bf CAT}(A, PB) \cong 
[ B^{\text op} \times A, {\bf Sets} ] = {\bf Prof}(A,B) \, .
\]
So the hom-categories can be identified with the familiar ones
in the definition of the bicategory $\bf Prof$, and 
then the Kleisli composition corresponds to
the standard composition in $\bf Prof$, and the
structure $2$-cells are as expected. So from Theorem 
\ref{kleislibicat} we deduce the following. 
\begin{proposition}
With structure as usually defined $\bf Prof$ is a bicategory.
\end{proposition}
This outline deduction is not as pointless as may appear. 
Where
in the literature can one find a proof?
The associativity may be made explicit,
but the reader is left to check the coherence
conditions. Another point is that the biequivalence
of $\bf Prof$ with the explicit sub-$2$-category of $\bf COC$ on
presheaf categories is usually given as extra information. 
For example it is proved
carefully in \cite{CW05}. However the moral of the treatment
here is that $\bf Prof$ is a bicategory in the first place 
just because it
is biequivalent to $\bf COC$.

\section{Distributivity}

Distributive laws between monads were introduced by Beck
who established the connection between lifting 
and composition of monads. The equally elementary connection
with extensions seems less well known but often arises in
semantics (see \cite{Hyland10}). From \cite{HNPR06}, I recall 
features of the general situation.
\begin{theorem}\label{kleisli}
Let $P$ and $S$ be monads on a category $\cal C$.
The following forms of data determine each other.
\begin {itemize}
\item A distributive law $SP \to PS$, that is, a natural transformation
preserving the units and multplications of
$P$ and $T$.
\item A lifting of $P$ along the forgetful functor $\Alg{S} \to {\cal C}$
to a monad $P^S$ on the category $\Alg{S}$.
\item An extension of $S$ along the free functor
${\cal C} \to \Kleisli{P}$ to a monad $S_P$ on the Kleisli category
$\Kleisli{P}$.
\end{itemize}
Moreover such data determines a monad structure on the
composite $PS$ with the category $\Alg{PS}$ of $PS$-algebras 
isomorphic to $\Alg{P^S}$ and the Kleisli
category for $\Kleisli{PS}$ isomorphic to $\Kleisli{S_P}$.
\end{theorem}
The subject of pseudo-distributive laws between
2-monads or more generally between pseudo-monads was already discussed
by Kelly~\cite{Kelly74} and has recently been the subject of
renewed attention, in \cite{Marm99} and \cite{Tanaka05}
for example. Work is required but unsurprisingly versions
of the basic categorical results go through.
Here I adapt the ideas to Kleisli structures to the extent needed
to arrive at Kleisli categories $\Kleisli{PS}$ of composed
Kleisli structures.

\subsection{Lifts of Kleisli structures}\label{lifts}
The situation which arises for us is the following.
We start with ${\cal A} \hookrightarrow {\cal K}$
a sub-$2$-category of a $2$-category. We have on the
one hand a Kleisli structure $P$ on ${\cal A} \hookrightarrow {\cal K}$
and on the other hand a 2-monad $(S, e, m)$ on $\cal K$ which restricts to a
2-monad on $\cal A$. (The less standard notation avoids confusion with
an earlier $\eta$.)

From $S$ we get a range of $2$-categories: there
are  $\Alg{S}_{\cal A}$ and  $\Alg{S}_{\cal K}$
the 2-categories of (strict) $S$-algebras
and pseudo-maps over $\cal A$ and $\cal B$; and
$\PsAlg{S}_{\cal A}$ and  $\PsAlg{S}_{\cal K}$,
the corresponding $2$-categories of pseudo-algebras.
The classic reference for material relating to
this kind of situation is \cite{BKP89}. In case $S$
is a flexible monad the choices here have little
importance, but our leading examples are not
flexible and so it makes best sense to
consider the inclusion
$\Alg{S}_{\cal A} \hookrightarrow \PsAlg{S}_{\cal K}$.
I display the obvious forgetful functors together
with the left pseudo-adjoints given by taking 
free $S$-algebras.
\begin{diagram}
 \Alg{S}_{\cal A}  & \rHook & \PsAlg{S}_{\cal K}\\
\uTo^{F} \dashv \dTo_{U} & & \uTo^{F} \dashv \dTo_{U} \\
{\cal A} & \rHook & {\cal K} & 
\end{diagram}
The details of all this are completely routine.

Now we consider {\em lifts} of the Kleisli structure
$P$ on ${\cal A} \hookrightarrow {\cal K}$ to a Kleisli
structure $P^S$on 
$\Alg{S}_{\cal A} \hookrightarrow \PsAlg{S}_{\cal K}$
By this we should presumably mean that the forgetful $2$-functors 
$U$ preserve the Kleisli structures up to coherent 
isomorphism. We can skip the details of the coherence as in all our examples,
$U$ will preserve the structure on the nose.

To give a lift of $P$ amounts in outline to the following.
For each $S$-algebra structure on $a$ we give an $S$-pseudoalgebra
structure on $Pa$ together with the structure of a pseudomap 
on $y_a: a \to Pa$;
and for each pseudomap with underlying $1$-cell $f:a \to Pb$
we give a pseudomap with underlying $1$-cell $f^\sharp : Pa \to Pb$,
in such a way that the $2$-cells $\eta_f$, $\kappa_A$ and $\kappa_{g,f}$
are pseudoalgebra $2$-cells. 

\subsection{Lifting presheaves}\label{liftingsec}
Our applications involve the presheaf Kleisli
structure $P$ on ${\bf Cat} \hookrightarrow {\bf CAT}$. To
give a lift of a $2$-monad $S$ is to give the following.
First for every $S$-algebra on a category $A$ we need to 
equip the presheaf category
$PA$ with the structure of an $S$-algebra in such a way
that the Yoneda preserves the structure in the up to coherent
isomorphism sense.
Secondly given $S$-algebras on categories $A$ and $B$  and an $S$-algebra 
pseudomap $f: A \to PB$ to the induced $S$-algebra $PB$, we
equip the left Kan extension $f^\sharp$ with the structure of
an $S$-algebra pseudomap so that we have an equality 
$f = y_A \cdot f^\sharp$ of pseudomaps.

It is as well to appreciate that none of this is automatic.
If $S$ gives structure which presheaves do not have then evidently
$P$ cannot lift. Biproducts gives an obvious example of that.
Also even if $S$ gives structure which presheaves do have
that structure may not be preserved by Yoneda.
Initial objects, coproducts, indeed most colimits are of this kind.
Finally even if $S$ satisfies the Yoneda condition, there may be a problem 
with the condition that appropriate left Kan extensions 
$f^{\sharp}: PA \to PB$ preserve $S$-structure. For example if $S$
is the monad which adds just equalizers or adds just pullbacks
then the left Kan extension condition fails.
In connection with the last point, I conjecture that there are very few
classes of finite limits which extend. A precise form of this is the
subject of a projected PhD thesis of Marie Bjerrum.

Fortunately there are many $2$-monads for which $P$ 
does lift. Examples include 2-monads for some familiar
classes of limits,
terminal object, products and finite limits.
The work of Im and Kelly \cite{IK86} provides additional examples: 
the $2$-monads for monoidal categories and for symmetric monoidal 
categories, together with minor variants. Finally there are
a range of miscellaneous examples: the $2$-monad equipping a category
with an endofunctor; the $2$-monad equipping a category with 
factorization. These last can be treated directly.
I expand a little on the examples which 
will be of most interest for the remainder of this paper.
\begin{enumerate}
\item The monad for a strict monoidal
category. Given a strict monoidal category $A$, $PA$ is
monoidal by the Day convolution tensor product
and so readily acquires
the structure of an $S$-pseudoalgebra. The Yoneda preserves
the tensor product in an up to coherent isomorphism sense,
and the left Kan extension condition is essentially in \cite{IK86}.
\item The case of the monad for a symmetric strict monoidal
category follows readily from the strict monoidal case, following
the lines of \cite{IK86}.
\item The monad for strictly associative
products could be another application
of \cite{IK86}, but can be handled directly.
Presheaf categories have products and
the Yoneda preserves them. Direct calculation
shows that if $f: A \to PB$ preserves
products then so does $f^{\sharp} : PA \to PB$.
(For $X \in PA$ the functor
$X \times - $ preserves colimits.)
\end{enumerate}

\subsection{Composed Kleisli structures}
I return to the general setting and suppose that $S$ 
is a $2$-monad and $P$ lifts from 
${\cal A} \hookrightarrow {\cal K}$ to $P^S$ on
$\Alg{S}_{\cal A} \hookrightarrow \PsAlg{S}_{\cal K}$. 
The aim of this section is to construct a 
Kleisli structure on ${\cal A} \hookrightarrow {\cal K}$
with the composite $PS$ as the basic operation on objects.

Applying our lifted Kleisli $P^S$ to free 
$S$-algebras $m_b : S^2 b \to Sb$
gives pseudoalgebras with $\sigma_b : SPSb \to PSb$ 
say as structure map.
In addition we have a pseudoalgebra map as in the diagram
\begin{diagram}
S^2B & \rTo^{Sy_{Sb}} & SPSb\\
\dTo^{m_b} & \cong & \dTo_{\sigma_b}\\ 
SB & \rTo_{y_{Sb}} & PSb
\end{diagram}

Now I define the unit structure $y^S$ and extension
structure $(-)^{\sharp^S}$ for a Kleisli
structure on $PS$. 
For each $a$ take as unit
\[
(y_a^S: a \to PSa) \quad = \quad (a 
\stackrel{e_a}{\longrightarrow} Sa 
\stackrel{y_{Sa}}{\longrightarrow} PSa) \, .
\]  
For $f : a \to PSb$ take as extension
\[
(f^{\sharp^S} :  PSa \to PSb )\quad = \quad (SA 
\stackrel{Sf}{\longrightarrow} SPSb
 \stackrel{\sigma_b}{\longrightarrow} PSb)^\sharp \, .
\]

Now for the structure $2$-cells. Given $f: a \to PSb$ 
let $\eta^S_f: f \to f^{\sharp^S}.y_a^S$ be the composite
of the $2$-cells isomorphisms
\[
f \cong \sigma_b . e_{PSb} . f \cong \sigma_b . Sf. e_a \cong
(\sigma_b . Sf)^\sharp. y_{Sa}. e_a \cong  f^{\sharp^S}.y_a^S
\]
To define $\kappa^S_a :(y_a^S)^{\sharp^S} \to 1_{PSa}$ 
take the composite of the $2$-call
isomorhisms
\[
(y_a^S)^{\sharp^S} \cong (\sigma_a . Sy_{Sa}. Se_a)^\sharp
\cong (y_{Sa} . m_a . Se_a)^\sharp \cong (y_{Sa})^\sharp \cong 1_{PSa} \, .
\] 
Finally given $f : a \to PSb$ and $g : b \to PSc$,  define
$\kappa^S_{g,f}: (g^{\sharp^S}.f)^{\sharp^S} \to g^{\sharp^S}.f^{\sharp^S}$
as the composite of the $2$-cell isomorphisms indicated.
\[
(g^{\sharp^S}.f)^{\sharp^S} \cong 
\big( \sigma_c . S(\sigma_c .Sg)^{\sharp^S} .Sf \big)^{\sharp^S}
\cong \big( (\sigma_c .Sg)^{\sharp^S} . \sigma_c .Sf \big)^{\sharp^S} \, ;
\]
\[
\big( (\sigma_c .Sg)^{\sharp^S} . \sigma_c .Sf \big)^{\sharp^S}
\cong (\sigma_c .Sg)^{\sharp^S}. (\sigma_c .Sg)^{\sharp^S} 
= g^{\sharp^S}. f^{\sharp^S} \, .
\]
In all cases the isomorphism in question is either structural,
comes from the Kleisli structure $P$ or is in the diagram above.
I did not make explicit what are the coherence isomorphisms
for a Kleisli structure and the reader will have to take on
trust my assertion that checking them is straightforward.
\begin{theorem}
Suppose that $S$ 
is a $2$-monad
and $P$ a Kleisli structure on ${\cal A} \hookrightarrow {\cal K}$
with $P$ lifting to $P^S$ on
$\Alg{S}_{\cal A} \hookrightarrow \PsAlg{S}_{\cal K}$.
Then $PS$ with structure above is a Kleisli structure on
${\cal A} \hookrightarrow {\cal K}$.
\end{theorem}
We do not need much more than this for applications so I say little
about $2$-dimensional versions of other points in Theorem \ref{kleisli}.
From a lift $P^S$ of $P$ we easily get a distributive law 
$\lambda : SP \to PS$
with components composites
\begin{diagram}
SPb & \rTo^{SP e_b} & SPSb & \rTo^{\sigma_b} & PSb & = & SPb & 
\rTo^{\lambda_b} & PSb
\end{diagram}
Preservation of structure by Yoneda specifies a natural isomorphism
\begin{diagram}
Sb & \rTo^{y_{Sb} } & PSb & \cong & Sb & \rTo^{Sy_b} & 
SPb & \rTo^{\lambda_b} & PSb 
\end{diagram}
which is part of the essential structure of a pseudo-distributivity.
I omit the rest of the structure.
Also there is an extension of $S$ to $S_P$ on $\Kleisli{P}$. In
principle it can be defined by composition with
$\lambda$ as in the ordinary case. A slicker equivalent approach
is to note that the forgetful $\Kleisli{U} : \Kleisli{P^S} \to \Kleisli{P}$
has a left biadjoint. The $S_P$ can be taken to be the
resulting pseudomonad on $\Kleisli{P}$. If we do this
carefully we get an isomorphism $\Kleisli{S_P} \cong \Kleisli{PS}$

\subsection{Extensions to Profunctors}
The bicategory $\bf Prof$ is very special and I give details of the
extension of $S$ on ${\bf Cat} \hookrightarrow {\bf CAT}$
to $S_P$ on $\bf Prof$. Extension is along 
$F: {\bf Cat} \to {\bf Prof}$ as in Section \ref{bicatsec}.
This pseudofunctor takes
$f:A \to B$ in $\bf Cat$ to $f_* : A \relto B$ in $\bf Prof$,
where $f_*(b,a) = B(b,fa)$. The first crucial
feature of the situation is that the arrow $f_*$ is a left adjoint
in $\bf Prof$ with right adjoint $f^*: B \relto A$ given by
$f^*(a,b) = B(fa, b)$.
Now extension gives for $f: A \to B$ a specified natural 
isomorphism between
\begin{diagram}
( Sf )_*  & = & SA & \rTo^{Sf} & SB & \rTo^{y_{SB} } &
 PSB
 \end{diagram}
 and 
 \begin{diagram}
 S_P(f_*) = SA & \rTo^{Sf} & SP & \rTo^{Sy_B} & SPB & \rTo^{\lambda_B} &PSB
 \end{diagram}
 Since we have specified adjunctions $(Sf)_* \dashv (Sf)^*$ and 
 $S(f_*) \dashv S(f^*)$, we get specified natural isomorphisms
 $(Sf)^* \cong S_P(f^*)$.   
 
In {\bf Prof}, there is a factorization of arrows analogous to
the factorization of relations through the graph of the relation
as the opposite of
a function followed by a function.
Given $M : A \relto B$ in $\bf Prof$, there is a category
$E(M)$ of elements of $M$ and functors $p: E(M) \to A$ and
$q: E(M) \to B$, which come together with an isomorphism
$M \cong q_*p^*$: in other words up to isomorphism
every arrow $M : A \relto B$ in $\bf Prof$ is a composite
$A \stackrel{p^*}{\relto} E(M) \stackrel{q_*}{\relto} B$.
Note that this means that the extension $S_P$ is determined
by the requirement that $S_P(f_*) \cong (Sf)_*$. For that 
requirement determines $S_P(f^*)$, and taken together
that determine the extension $S_P$.

\section{Algebraic theories}

For the remainder of this paper I shall work in the setting provided
by the presheaf Kleisli structure $P$
on ${\bf Cat} \hookrightarrow {\bf CAT}$. I consider notions of
algebraic theory determined by $2$-monads $S$ for which we
have a lift $P^S$ of $P$. The action takes place in the
Kleisli bicategory $\Kleisli{PS}$. I shall write its maps
as $A \relto SB$ which makes visible the identification
with $\Kleisli{S_P}$. (At a few points I shall make active use
of the extension $S_P$ of $S$.)

Since in many modern applications one needs enriched notions of algebraic
theory I note that the theory described enriches
readily. The background needed for that is contained in \cite{Kelly82}.
There are details to check but no essential difficulties appear.
That of course is to say that the basic enriched theory is already
apparent in what I present.

\subsection{Kleisli objects in profunctors}\label{kleisliobject}
My approach to algebraic theories uses monads in a
Kleisli bicategory $\Kleisli{PS}$
so to prepare the ground I say a little about monads in $\bf Prof$. A monad 
is given by a profunctor $M: A \relto A$ with a unit $2$-cell
$I_A \cell M$ and composition $2$-cell $M\cdot M \cell M$
satisying the usual equations. For $\bf Prof$ the Eilenberg-Moore
and Kleisli objects for a monad coincide, but it is the Kleisli
which is important here. The universal diagram is of
the form
\begin{diagram}
A & & \rRel^M & & A \\
 & \rdRel_{k_*} & \Leftarrow & \ldRel_{k_*} & \\
 & & C(M) & & 
\end{diagram} 
where the $2$-cell satisfies the conditions for
the algebra for a monad. Concretely the Kleisli object 
$C(M)$ has the same objects as $A$ and its arrows
$C(M)(b,a) = M(b,a)$ are given by the elements of $M$.
Composition and units come from the corresponding structure
on $M$, and the comparison arrow $k_* : A \to C(M)$
corresponds to an identity on objects functor $k : A \to C(M)$.

Every functor can be factorised as an identity on objects
followed by a full and faithful functor. This well known
factorisation system on $\bf Cat$ can be derived from
the Kleisli construction. For let $f : A \to B$ 
be a functor. As $f_* \dashv f^*$, the composite 
$f^* \cdot f_* : A \relto A$ is a monad. Moreover the
counit of the adjunction gives a $2$-cell
\begin{diagram}
A & & \rRel^{f^* \cdot f_*} & & A \\
 & \rdRel_{f_*} & \Leftarrow & \ldRel_{f_*} & \\
 & & B & & 
\end{diagram}
satisfying the algebra conditions. So by the
universality of the Kleisli construction, 
$A \stackrel{f}{\longrightarrow} B$ factors uniquely
as $A \stackrel{f}{\longrightarrow} 
C(f^* \cdot f_*) \stackrel{}{\longrightarrow} B$.

\subsection{Intuition on composed Kleisli bicategories}
Let $S$ be a $2$ monad with $P$ lifting so that we have the
composed Kleisli bicategory
$\Kleisli{PS} \cong \Kleisli{S_P}$ with arrows $A \relto SB$. 
I shall use $\opcomm$ to denote
the composition in $\Kleisli{PS}$: that is more or less in accord 
with practice in the operads community.

Let me give an intuitive syntactic 
reading of the arrows and their composition
in $\Kleisli{PS}$. Consider $F : A \relto SB$. Write $a \in A$ for 
objects of $A$
and ${\bf b} \in SB$ for objects of $SB$. Then an
element $f \in M({\bf b}, a)$ should be thought
of as a formal function or function symbol
with input arity $\bf b$ and
output arity $a$. (The choice of $S$ determines the nature
of the input arities.) Now suppose we have $F : A \relto SB$ and
$G : B \relto SC$. By definition $G \opcomm F: A \relto SC$ is 
the composite $A \stackrel{F}{\longrightarrow} 
SB \stackrel{SG}{\longrightarrow} S^2C \stackrel{m_*}{\longrightarrow} SC$.
What is an element of  $G \opcomm F ({\bf c}, a)$? In the language 
of function symbols the intuition is that it is determined
by a function symbol $f \in F({\bf b}, a)$ together with a certain $S$-arities
worth of functions symbols $g_i \in G({\bf c}_i, b_i)$, with the $b_i$
making up $\bf b$; the $g_i$ should be thought of as substituted
into the input arity $\bf b$ and the arities ${\bf c}_i$ composed together
by the multiplication in $S$ to give the single arity ${\bf c}$.
The overall idea is that the $2$-monad $S$ has determined
a general notion of substitution.

\subsection{Theories as monads}
In this section I explain the idea that
a many-sorted $S$-algebraic theory or 
simply an $S$-multicategory,
can be represented as a monad 
$M: A \relto SA$ in $\Kleisli{PS}$, or equivalently as monoids in
the endocategory $\Kleisli{PS}(A,A)$. The familiar single-sorted 
theories correspond to the case $A=1$.
I shall use the term monad here though the
operads community often says monoid.

A monad in $\Kleisli{PS}$ consists of a 
profunctor $M: A \relto SA$ with a unit 
and composition $2$-cells
$\eta_{A*} \cell M$ and $M \opcomm M \cell M$
satisfying the usual equations.
In terms of the syntactic reading, the identity provides identity
functions of input arity $e(a)$ and output arity $a$. This is usually
represented in syntax by a variable. Composition gives the interpretation of 
formal composites of $f \in M({\bf a}, a)$ with an $S$-indexed family
${\bf g} \in SM(\alpha, {\bf a})$,
the interpretation being given in $M(m(\alpha),a)$. So it provides 
the interpretation of composite symbols. That is in general
terms exactly what we expect of an algebraic theory. I shall tighten
up this idea shortly but for the moment let us assume that an
$S$-algebraic theory is just a monad in $\Kleisli{PS}$

To define a map between $S$-algebraic theories we need to exploit
the pseudofunctor $F: {\bf Cat} \to {\bf Prof}$: we cannot simply
take monad maps in $\Kleisli{PS}$ as they evidently do not give
what we want. 
\begin{definition}
A {\em map of theories} from $M: A \relto SA$ to $N: B \relto SB$
consists of a functor $f: A \to B$ together with a $2$-cell
\begin{diagram}
A & \rRel^{M} & SA \\
\dRel^{f_*} & \Downarrow & \dRel_{Sf_*}\\
B & \rRel_{N} & SB
\end{diagram}
which satisfies commutation conditions expressing compatibility
with the unit and composition $2$-cells. A {\em $2$-cell between maps} 
is given by a natural transformation
between the functors compatible with the the $2$-cells in
the obvious sense. Taken together this data gives us a $2$-category
$\Mult{S}$ of {\em $S$-algebraic theories}.
\end{definition}
I need immediately to refine that definition.
It is almost right except for the fact 
that $A \stackrel{M}{\relto} SA$ carries not just the data of the
multicategory but data corresponding to the category $A$ and its
action on the multicategory. This data is unwelcome
and there are two ways to suppress it. 
\begin{enumerate}
\item One can either insist that the categories 
$A$ are discrete so that there is no information either in $A$ or its action;
This is the evident choice.
\item One can require that the category $A$ is exactly the underlying
category of the multicategory $M$, i.e. it is {\em normal}
in Australian terminology.
\end{enumerate} 
For the theory developed here the second is the
more natural choice so let us assume normal in the definition
above. For clarity I explain precisely what it means
for a monoid $A \stackrel{M}{\relto} SA$ to be normal. The unit $2$-cell
$\eta_{A*} \cell M$ gives by transpose a cell $I_A \cell \eta_A^*\cdot M$,
and we require that this be an isomorphism. The discussion 
in Section \ref{kleisliobject} shows that this what is intended. 
The point is that $\eta_A^* \cdot M$ becomes a monad in $\bf Prof$ 
and so gives an identity on objects functor from $A$ to what is 
the underlying category of $M$; and we want this to be an isomorphism.

I do not want to make a big deal about the choice here. The rest of 
the paper can be read without much change with the other choice.
In any case the existence of two distinct embeddings, one discrete,
the other chaotic, is a common situation
to which Lawvere has frequently drawn attention.

\subsection{Models of algebraic theories}

This is a rich topic and I just give the basic definition
and make a few remarks about the monadic approach to algebra.
First I recall an argument from \cite{HP07} in favour of the
Lawvere theory approach to algebraic theories over the monadic. 
The monadic approach gives immediate access only to what are
called algebras for the monad, that is, models in
the ambient category of the monad. The Lawvere theory
approach from \cite{Lawvere63} gives models in any category with products.
In line with the general practice of categorical logic,
models are product preserving functors. How does the
bicategorical approach to
algebraic theories compare?

The simple view is that models of an $S$-algebraic
theory $A \stackrel{M}{\to} SA$ can be taken in any 
$S$-algebraic theory. A category of models in 
$B \stackrel{N}{\to} SB$ is just
the category of maps of algebraic theories from $A \stackrel{M}{\to} SA$
to $B \stackrel{N}{\to} SB$. That is the fundamental notion, but
perhaps it will leave readers unsatisfied. What about models in
categories, or rather, taking the lead from \cite{Lawvere63}, in categories
with $S$-structure up to isomorphism, that is to say in 
$S$-pseudoalgebras? Happily there is an elegant account of that.

First suppose that we are given an $S$-algebra $a: SA \to A$. This gives rise
to an $S$-multicategory in a ridiculously simple way. One 
takes the profunctor $A \stackrel{a^*}{\relto} SA$ and equips it with a monoid
structure. The unit $\eta_* \cell a^*$ is the transpose of the
isomorphism (identity) $a_* \eta_* \cell 1_A$.
The multiplication $a^* \opcomm a^* \cell a^*$, that is,
$\mu_* (Sa^*) a ^* \cell a^*$ is the transpose of the isomorphism
$(Sa^*) a ^* \cell \mu^* a^*$. It is straightforward to
check that this makes $A \stackrel{a^*}{\relto} SA$ 
a normal $S$-multicategory.
This construction extends to a $2$-functor $\Alg{S}$ to $\Mult{S}$.
We shall consider it further in the next section.

However nothing depended on the supposition that $a: SA \to A$ was a strict 
as opposed to a pseudoalgebra; and nothing beyond a question of
extending definitions on the supposition that $A$ was small.
Certainly for any $x : SX \to X$ in $\PsAlg{S}_{\bf CAT}$,
we can formally construct a large $S$-multicategory
$X \stackrel{x^*}{\to} SX$. Then by definition a model of
$S$-algebraic theory $A \stackrel{M}{\to} SA$ in an $S$-pseudoalgebra
$x : SX \to X$ is just a map of theories $A \stackrel{M}{\to} SA$
to $X \stackrel{x^*}{\to} SX$, and we get a (possibly large)
category of models.

I close this section with observations from the monadic
point of view. It is clear that a monad in a bicategory acts
by composition (on either side) as a monad on suitable
hom-categories. In the case of our leading examples, this
gives an immediate explanation of the monad on $\bf Sets$
generated by a single-sorted algebraic theory $1 \stackrel{M}{\to} S1$.
In $\Kleisli{S_P}$, $M$ is a monad on $1$, and so
we obviously have an action by composition 
in particular on $\Kleisli{S_P}(1,0)$. But
in each case we have $S0 \cong 1$, so that 
$\Kleisli{S_P}(1,0)$ is isomorphic to $\bf Sets$. Thus we get 
the usual monad
on $\bf Sets$.

From the bicategorical point of view the choice of $\Kleisli{S_P}(1,0)$
is rather arbitrary. There is at least one other compelling
case which is $\Kleisli{S_P}(1,1)$. I note that when $S$ is 
symmetric strict monoidal categories $2$-monad, this choice gives
the notion of twisted (French tordue) algebra.
Clearly there is more to be said to reconcile the different points
of view on models.

\subsection{The free $S$-algebra}\label{freesec}
In the previous section we saw the simple functor $\Alg{S} \to \Mult{S}$
giving the $S$-algebraic theory corresponding to an $S$-algebra.
We expect a left adjoint, that is, 
given an $S$-multicategory we expect to be able to construct freely from it
an $S$-algebra. We can do that very generally and I sketch the
essential point. Suppose that we have a multicategory
$M : A \relto SA$ and $S$-algebra $b: SB \to B$. Then there is 
a correspondence between diagrams of the two forms
\begin{diagram}
A & \rRel^M & SA & \qquad \qquad & SA & \rRel^{S_PM} & S^2A &
\rRel^{(m_A)_*} & SA \\
\dRel^{f_*} & \Downarrow & \dRel_{Sf_*} &&& \rdRel_{g_*}& \Leftarrow &
\ldRel_{g_*} & \\
B & \rRel_{b^*} & SB & & & & B & & 
\end{diagram}
where predictably $g:SA \to B$ and $f: A \to B$ determine each
other by $g=b.Sf$ and $f = g.e_A$. In the diagram on the left
we take the structure of a map of multicategories from 
$M : A \relto SA$ to $b^* : B \relto SB$, while on the left
we take a Kleisli cone in a prima facie non-evident
bicategory {\em of lax $S$-profunctors} in which objects
are $S$-algebras and arrows profunctors preserving
$S$-algebra structure in a lax sense. (This precise
formulation is Richard Garner's.) It follows that quite generally
the free $S$-algebra on an $S$-multicategory is given by a special
Kleisli construction. 
\begin{theorem}\label{freeSalg} If Kleisli 
objects exist in the bicategory of lax
$S$-profunctors, then there is a left adjoint $\Mult{S} \to \Alg{S}$ to
$\Alg{S} \to \Mult{S}$.
\end{theorem}
As with many colimits the existence of Kleisli
objects is shown by an iterative construction. (In view of
Proposition \ref{comonad} it is not surprising that this is very close
to a construction in \cite{BKP89}.)
However in our leading examples there is a much simpler approach. 
All satisfy the condition that the $2$-monad $S$ 
preserves bijective on objects functors.
It is easy to check the following basic observation.
\begin{proposition}
Let $S$ be a $2$-monad on ${\bf Cat} \hookrightarrow {\bf CAT}$
with extension $S_P$ on $\bf Prof$. Then $S$ preserves
bijective on objects functors if and only if $S_P$ preserves Kleisli objects. 
\end{proposition}

Now suppose we are in the situation above.
Let $M: A \relto SA$ be an $S$-multicategory and consider
the composite 
\begin{diagram}
SA & \rRel^{SM} & S^2 A & \rRel^{\mu_*} & SA
\end{diagram}
which is certainly a monad in $\rm Prof$. We can construct its Kleisli
object $C(M)$ as in the diagram
\begin{diagram}
 SA & \rRel^{S_PM} & S^2A & \rRel^{(m_A)_*} & SA \\
& \rdRel_{h_*}& \Leftarrow & \ldRel_{h_*} & \\
 & & C(M) & & 
\end{diagram}
and exploit the basic fact that we have maps where indicated.
Applying $S_P$ to this diagram we get another Kleisli
object and so by the universality explained earlier 
a map $SC(M) \to C(M)$. Easily that makes $C(M)$ into an
$S$-algebra and $h : SA \to C(M)$ a strict map of $S$-algebras.
It is now routine to check the following.

\begin{proposition}
Suppose $S$ is a $2$-monad on ${\bf Cat} \hookrightarrow {\bf CAT}$
extending to $S_P$ on $\rm Prof$, and suppose further that $S$ preserves
bijective on objects functors. Then for every monad 
$M: A \relto SA$ in $\Kleisli{PS}$,
its Kleisli object is given by $C(M)$ with the structure above.
\end{proposition}

\subsection{Lawvere theories and PROPs}

I record what the construction of the free $S$-algebra amounts
to in our leading special examples from Section \ref{liftingsec}, 
which I now take in reverse
order.

First take $S$ to be the monad for strict products.
An $S$-algebraic theory $1 \stackrel{M}{\relto} S1$ corresponds exactly
to what is usually called an algebraic theory: my formulation is a fancy way to
describe the notion of abstract clone from universal algebra.
If $M$ is an algebraic theory then the construction of $C(M)$ 
gives exactly the corresponding Lawvere theory.
More general $S$-algebraic theories $A \stackrel{M}{\relto} SA$
might naturally be called many-sorted algebraic theories.
I prefer and shall use the precise terminology cartesian
theories. For $M$ a cartesian theory $C(M)$ gives a 
category with strict products.
By Section \ref{freesec} the models of $C(M)$ in the usual 
sense of categorical logic are
equivalent to the models of the cartesian theory.

Now take $S$ to be the monad for
symmetric strict monoidal categories.
An $S$-algebraic theory $1 \stackrel{M}{\relto} S1$ corresponds exactly
to a symmetric operad. In this case the construction of $C(M)$
gives the corresponding
$PROP$, in effect symmetric monoidal category generated
by one object. A more general $S$-algebraic theory 
$A \stackrel{M}{\relto} SA$ 
is what is called a coloured operad.
For those $C(M)$ gives a symmetric strict monoidal
category whose models in the sense of symmetric monoidal categories
are models for the coloured operad.

Finally take $S$ to be the monad for
strict monoidal categories.
An $S$-algebraic theory $1 \stackrel{M}{\relto} S1$ corresponds exactly
to a non-symmetric operad. In this case the construction of $C(M)$
gives the corresponding non-symmetric
$PROP$. (That is a slight contradiction of terminology!)
A more general $S$-algebraic theory $A \stackrel{M}{\relto} SA$ 
is a non-symmetric coloured operad.
For those $C(M)$ gives a strict monoidal
category whose models in the sense of monoidal categories
are models for the coloured non-symmetric operad.

All the above is quite straightforward, but
there is a significant difference between the
first case and the other two.
This is most easily explained by making precise a
connection with \cite{BKP89}. For a $2$-monad $S$,
one has $2$-categories, $\Alg{S}_{str}$, $\Alg{S}$ and $\Alg{S}_{lax}$
of strict $S$-algebras with strict, pseudo and lax maps respectively.
One has forgetful functors
$\Alg{S}_{str} \to \Alg{S} \to \Alg{S}_{lax}$.
A major achievement of \cite{BKP89} is the construction of
left adjoints
\[
(-)' : \Alg{S} \to \Alg{S}_{str} \qquad {\rm and} \qquad
(-)^{\dagger} : \Alg{S}_{lax} \to \Alg{S}_{str} \, .
\]
\begin{proposition}\label{comonad}
Suppose $S$ is a $2$-monad with extension $S_P$. Then for every
$S$-algebra $a: SA \to A$,  the free $S$-algebra generated by 
the $S$-algebraic theory
$A \stackrel{a^*}{\to} SA$ gives the left adjoint $(A)^{\dagger}$.
\end{proposition}
Another way to say the same thing is that the $2$-comonad on
$\Alg{S}_{str}$ generated by the adjunction of Theorem \ref{freeSalg}
coincides with that given in \cite{BKP89}.
The following is then immediate.
\begin{proposition}
The $2$-category of $\Mult{S}$ of $S$-algebraic theories
embeds in the $2$-category of $(-)^{\dagger}$-coalgebras.
\end{proposition}

Now suppose in our standard setting that the $2$-monad $S$ 
is colax idempotent. Then lax $S$-algebra maps are automatically 
pseudo and we can identify the $2$-functors
$(-)^{\dagger}$ and $(-)'$. Now by \cite{BKP89} each $S$-algebra $A$
is equivalent to $A'$ which is a free $(-)'$-coalgebra, equivalently
a free $(-)^{\dagger}$-coalgebra. But the category of 
free $(-)^{\dagger}$-coalgebra is just $\Alg{S}_{lax}$.
Hence every $S$-algebraic theory is equivalent to one
arising from an $S$-algebra.
Since the $2$-monad for strict products is colax idempotent,
it follows that up to equivalence Lawvere theories
are just algebraic theories. In the many object case, categories with
products can be identified with cartesian theories.

For general $S$ however the above analysis does not apply.
Specifically, the $2$-monad for symmetric
strict monoidal categories is certainly not colax idempotent
and it is not the case that every symmetric strict monoidal
category arises up to equivalence,
as the free symmetric monoidal category generated by
a coloured operad. For example
the PROP for comonoids does not so arise. Why? Because
if a theory is given by an operad, then there is a corresponding monad
on any symmetric monoidal closed cocomplete category and there is no free 
comonoid functor on most such categories. In particular there is none
in the case of ${\bf Vect}_k$, the category of vector spaces over
a field $k$ with its standard tensor structure. Why? Because the
forgetful does not preserve limits. Why? Well the terminal coalgebra
is just the field $k$ with trivial coalgebra structure and $k$ is
not terminal in ${\bf Vect}_k$.

\section{Comparing notions of theory}

The aim of this section is to develop the theory needed
to support the free and forgetful functors linking
various notions of theory. Since notions of theory here
correspond to certain $2$-monads, this depends on understanding
maps between them.

\subsection{Compatible maps of monads}
Suppose now that $S$ and $T$ are $2$-monads on 
${\bf Cat} \hookrightarrow {\bf CAT}$
with liftings $P^S$ and $P^T$ of the presheaf
Kleisli structure to the categories of algebras. 
Suppose $k:T \to S$ is a map of $2$-monads,
inducing maps $\Alg{S}_{\bf Cat} \to \Alg{T}_{\bf Cat}$ and
$\PsAlg{S}_{\bf CAT} \to \PsAlg{T}_{\bf CAT}$, which is in addition compatible
with taking categories of presheaves in that we have
a natural isomorphism
\begin{diagram}
\Alg{S}_{\bf Cat} & \rTo^P & \PsAlg{S}_{\bf CAT} \\
\dTo & \cong & \dTo \\
\Alg{T} _{\bf Cat}& \rTo_P & \PsAlg{T}_{\bf CAT} \, .
\end{diagram}
We have the distributivities
$\lambda^S: SP \to PS$ and $\lambda^T : TP \to PT$. Evaluating the isomorphism
at free $S$-algebras gives isomorphisms between the resulting families
of $T$-algebras
\begin{diagram}
TPSA & \rTo^{k_{PSA}} & SPSA & \rTo^{\lambda^S_{SA}} &
PS^2A & \rTo^{Pm^S} & PSA \\
TPSA & \rTo^{\lambda^T_{SA}} & PTSA & \rTo^{Pk_{SA}} &
PS^2A & \rTo^{Pm^S} & PSA \\
\end{diagram} 
with underlying category $PS$. That gives a choice of natural isomorphism
\begin{diagram}
TP & \rTo^{kP} & SP \\
\dTo^{\lambda^T} & \cong & \dTo_{\lambda^S} \\
PT & \rTo_{Pk} & PS
\end{diagram}
It is relatively straightforward to see that conversely this 
condition implies the compatibility of the maps 
above with the taking of presheaves.

It will be good at this point to see that we can find simple
cases of compatible maps of monads between our leading
examples.

Consider first the case when $S$ and $T$ are the $2$-monads 
for categories with strict finite products and strict
symmetric monoidal categories respectively. Any category with 
strict products is symmetric strict monoidal,
so we certainly have a map $k:T \to S$ of
$2$-monads. What about the compatibility with presheaves?
Well let $A$ be a category with strict finite products.
Going one way round the natural isomorphism square we get
first $PA$ as a category with finite products (considered
as an $S$-pseudoalgebra) and then $PA$ again and with the
products structure simply considered as (symmetric) monoidal
structure (as a $T$-pseudoalgebra).
Going the other way round we get first $A$ considered
now as a symmetric strict monoidal category, and then
$PA$ with the induced Day tensor product. Now it is both folklore
and an easy computation to show that the Day tensor product
obtained from cartesian product structure is itself
a cartesian product. So there is an isomorphism
given by the unique isomorphism between choices of
products in the compatibility square. The uniqueness
makes naturality rather evident.

Now what about the case when $S$ and $T$ are the $2$-monads for
symmetric strict monoidal categories and strict monoidal categories
respectively. Here the situation is even easer. One way round
we take the Day tensor product and forget the symmetry;
the other we forget the symmetry and then take the Day tensor
product. That amounts to exactly the same thing, so quite
exceptionally the natural isomorphism is the identity.
(Note the first case could have been made similarly easy by using the
Day tensor product throughout, but that is rather unilluminating.)

\subsection{Comparison of Kleisli bicategories}

At this point it is probably simplest to think in terms of the
extended pseudomonads $T_P$ and $S_P$.
Take $f: A \relto B$ in {\bf Prof}, that is, $f : A \to PB$ in $\bf CAT$.
The images of $f$ under $T_P$ and $S_P$
lie in a diagram
\begin{diagram}
TA & \rTo^{Tf} & TPB & \rTo^{\lambda^T_B} & PTB \\
\dTo^{k_A} & \cong & \dTo^{k_{PB}} & \cong & \dTo_{Pk_B} \\
SA & \rTo_{Sf} & SPB & \rTo_{\lambda^S_B} & PSB 
\end{diagram}
Rotating for readability that gives isomorphisms
\begin{diagram}
T_PA & \rRel^{k_*} & S_PA \\
\dRel^{T_Pf} & \cong & \dRel_{S_Pf}\\
T_PB & \rRel_{k_*} & S_PB 
\end{diagram}
which gives the family $(k_A)_* : TA \relto SA$ the structure
of a strong transformation $k_* : T_P \to S_P$ on 
$\bf Prof$. It follows automatically, or if you prefer
it can be proved 
directly, that
the family of right adjoints $(k_A)^* : SA \relto TA$ have the structure
of a lax transformation $k^* : S_P \to T_P$ on 
$\bf Prof$. The coherent $2$-cells for the lax transformation
are as in the diagram
\begin{diagram}
S_PA & \rRel^{k^*} & T_PA \\
\dRel^{S_Pf} & \Downarrow & \dRel_{T_Pf}\\
S_PB & \rRel_{k^*} & T_PB 
\end{diagram}

Now consider composition with $k_*$ and $k^*$ as operations
on the Kleisli bicategories. The former
\begin{diagram}
\hat{k}_* \, : \, (A \stackrel{u}{\relto} TB) & \qquad \longmapsto \qquad & 
(A \stackrel{u}{\relto} TB \stackrel{k_*}{\relto} SB )
\end{diagram}
gives a pseudofunctor $\hat{k}_* : \Kleisli{PT} \to \Kleisli{PS}$,
while the latter
\begin{diagram}
\hat{k}^* \, : \, (A \stackrel{v}{\relto} SB ) & \qquad \longmapsto \qquad & 
(A \stackrel{v}{\relto} SB\stackrel{k^*}{\relto} TB )
\end{diagram}
gives a lax functor $\hat{k}^*: \Kleisli{PS} \to \Kleisli{PT}$.
The adjunction $k_* \dashv k^*$ induces isomorphisms
between $\Kleisli{S_P}(k_* u,v)$ and $\Kleisli{T_P}( u, k^*v)$
so $\hat{k}^*$ is locally right adjoint to $\hat{k}_*$ in
the sense that we have adjunctions
\[
\xymatrix{
   {\Kleisli{PT}}(A,B) 
 \ar@<1.2ex>[rr]^(.52){\hat{k}_*}
 \ar@{}[rr]|{\bot}    
 & & \ar@<1.2ex>[ll]^{\hat{k}^*} {\Kleisli{SPS}}(A, B) } \, .
\]

\subsection{Comparison of notions of theory}

The simple fact that lax functors between bicategories
(morphisms in Benabou's terminology) take monads to monads
already plays an important role in the original paper
\cite{Benabou67}. As the condition of being normal
is algebraic it is easy to see that both $\hat{k}_*$ and
$\hat{k}^*$ preserve algebraic theories. Maps require more
care: despite the importance of lax functors to the Benabou
vision, they do not preserve maps of monads. (In a sense that is part of
the point of \cite{Benabou67}.) Also one has to pause
because maps of algebraic theories correspond to rather
special maps of monad, using the free pseudo-functor
${\bf Cat} \to \Kleisli{PS}$. So there is a little work to do.

The case of the left adjoint $\hat{k}_*$ is straightforward
as it is a pseudo-functor, and its action on maps in
$\Kleisli{PS}$ coming from $\bf Cat$ is essentially
trivial. So the action of $\hat{k}_*$
\begin{diagram}
A & \rRel^{M_T} & TA & & A & \rRel^{M_T} & TA & \rRel^{k_{A*}} & SA \\
\dRel^{f_*} & \Downarrow & \dRel_{Tf_*} & \qquad \mapsto \quad &
 \dRel^{f_*} & \Downarrow & \dRel^{Tf_*} & \cong & \dRel_{Sf_*} \\
B & \rRel_{N_T} & TB & & B & \rRel_{N_T} & TB & \rRel_{k_{B*}} & SB
\end{diagram}
takes maps of $T$-algebraic theories to maps of $S$-algebraic theories.
That action extends to $2$-cells so we get a $2$-functor
$\hat{k}_* : \Mult{T} \to \Mult{S}$.

Prima facie it looks as if we need more work to deal with
the lax functor $\hat{k}^*$, but that is counter intuitive
as the local right adjoint $\hat{k}^*$ is foregtful and
so should be straightforward. That intuition is correct
and one can argue concretely since
for $A \stackrel{M}{\relto} SA$, we have 
$\hat{k}(M) ({\bf a}, a) \cong M(k{\bf a}, a)$. Naturally
there is a more abstract apporach but either way
we get the following.
\begin{theorem}
Composition with $k_*$ and $k^*$ induce a $2$-adjunction
\[
\xymatrix{
   {\Mult{T}} 
 \ar@<1.2ex>[rr]^(.52){\hat{k}_*}
 \ar@{}[rr]|{\bot}    
 & & \ar@<1.2ex>[ll]^{\hat{k}^*} {\Mult{S}} } \, .
\]
\end{theorem}
Here the right adjoint is the easy forgetful $2$-functor,
while the left adjoint is the more subtle free $2$-functor.
One gets exactly what is expected in the examples.
One can forget from an algebraic theory to a symmetric operad,
and there is a free algebraic theory generated by such an
operad. Similarly one can forget from a symmetric operad
to a non-symmetric operad and there is a free symmetric operad
on a non-symmetric one. (Of course there are the composites.)

I close with a remark from Tom Leinster about a difference between the
two cases. In the first the free functor $\hat{k}^*$ reflects isomorphisms.
If a cartesian theory arises from a coloured operad then the 
operad in question is determined up to isomorphism. Restricting
to a single sort, being (symmetric) operadic is thus a property of an algebraic
theory. However in the second case the free 
functor $\hat{k}^*$ does not reflect isomorphisms. An example
of non-isomorphic non-symmetric operads which generate
isomorphic symmetric operads
is given in \cite{Leinster06}. Hence there is an issue if one says 
of a symmetric operad that it is non-symmetric. Probably one
should make a choice of the non-symmetric operad in question,
in which case being non-symmetric is not a property but
additional structure.

\end{document}